\documentclass[12pt]{amsart}

\makeatletter

\providecommand{\LyX}{L\kern-.1667em\lower.25em\hbox{Y}\kern-.125emX\@}

 \theoremstyle{plain}    
 \newtheorem{thm}{Theorem}[section]
 \numberwithin{equation}{section} 
 \numberwithin{figure}{section} 
 \theoremstyle{plain}    
 \newtheorem{cor}[thm]{Corollary} 
 \theoremstyle{plain}    
 \newtheorem{lem}[thm]{Lemma} 
 \theoremstyle{plain}    
 \newtheorem{prop}[thm]{Proposition} 
 \theoremstyle{definition}
 \newtheorem{defn}[thm]{Definition}
 \theoremstyle{definition}
  \newtheorem{example}[thm]{Example}
 \theoremstyle{definition}
  \newtheorem{condition}[thm]{Condition}
 \theoremstyle{remark}    
 \newtheorem*{acknowledgement*}{Acknowledgement} 

\newcommand{\Tr}{\operatorname{Tr}}
\newcommand{\dist}{\operatorname{dist}}
\newcommand{\id}{\operatorname{id}}
\makeatother
\begin{document}
\title[Some estimates for $\delta^*$]{Some Estimates for non-microstates free entropy dimension, with applications to $q$-semicircular families}

\address{Department of Mathematics, University of California, Los Angeles, CA 90095, USA}

\email{shlyakht@math.ucla.edu}

\author{Dimitri Shlyakhtenko}

\begin{abstract}
We give an general estimate for the non-microstates free entropy dimension
$\delta ^{*}(X_{1},\ldots ,X_{n})$. If $X_{1},\ldots ,X_{n}$ generate
a diffuse von Neumann algebra, we prove that $\delta ^{*}(X_{1},\ldots ,X_{n})\geq 1$.
In the case that $X_{1},\ldots ,X_{n}$ are $q$-semicircular variables
as introduced by Bozejko and Speicher and $q^{2}n<1$, we show that
$\delta ^{*}(X_{1},\ldots ,X_{n})>1$. We also show that for $|q|<\sqrt{2}-1$,
the von Neumann algebras generated by a finite family of $q$-Gaussian
random variables satisfy a condition of Ozawa and are therefore solid:
the relative commutant of any diffuse subalgebra must be hyperfinite.
In particular, when these algebras are factors, they are prime and
do not have property $\Gamma $.
\end{abstract}
\maketitle

\section{Introduction.}

In \cite{Speicher:q-comm}, Bozejko and Speicher introduced a deformation
of a free semicircular family of Voiculescu \cite{DVV:free,DVV:circular},
parameterized by a number $q\in [-1,1]$. Their $q$-semicircular
family $X_{1},\ldots ,X_{n}$ is represented on a deformed Fock space
and generates a finite von Neumann algebra, which is non-hyperfinite
for $n\geq 2$ and $q\in (-1,1)$ \cite{nou:qNonInjective}. For $q=0$,
$X_{1},\ldots ,X_{n}$ are semicircular and free in the sense of Voiculescu's
free probability theory. For $q=1$, they form a family of independent
Gaussian random variables, and for $q=-1$ they generate a Clifford
algebra related to the Canonical Anti-Commutation Relations.

The von Neumann algebra $M_{q}=W^{*}(X_{1},\ldots ,X_{n})$ generated
by such a $q$-semicircular family remains a mystery. One would like
to decide whether $M_{q}$ falls into the world of {}``free'' von
Neumann algebras (as $M_{0}\cong L(\mathbb{F}_{n})$ does), or whether
$M_{q}$ move away from that world immediately once $q\neq 0$. 

In this paper, we give evidence for the fact that $M_{q}$ stay close
to free group factors. We show that Voiculecsu's non-microstates free
entropy dimension computed on the $q$-semicircular generators of
$M$ is $>1$ for small $q$. Furthermore, we show that $M_{q}$ for
small $q$ satisfy the Ozawa condition \cite{ozawa:solid}, implying
that they are {}``solid'', in particular, prime. We still don't
know if they are always factors, except when this is guaranteed by
the results of Sniady \cite{sniady:qfactor}.

It is curious that the essential ingredient in both the computation
of free entropy dimension and in the proof of Ozawa's condition is
the existence of certain operators $r_{1},\ldots ,r_{n}$ satisfying
$[r_{j},X_{i}]=\delta _{ij}\Xi _{j}$, with $\Xi _{j}$ compact or
even Hilbert-Schmidt. The estimate on free entropy dimension follows
roughly Voiculescu's argument that the existence of a {}``dual system''
for $Y_{1},\ldots ,Y_{n}$ implies that their free entropy dimension
is large.

Our proof of Ozawa's condition gives another argument for why it is
satisfied by the free group von Neumann algebra; our argument, strangely
enough, avoids all use of the free group or of its boundary. It would
be very interesting to understand if our proof gives another point
of view on the action of the free group on its boundary, which is
more amenable to generalizations.

We also give a general lower bound on Voiculescu's non-microstates
free entropy dimension of a self-adjoint $n$-tuple $(X_{1},\ldots ,X_{n})$
in terms of the Murray-von Neumann dimension of the space of certain
derivations from the algebra generated by $X_{1},\ldots ,X_{n}$ into
Hilbert-Schmidt operators. As an application, we prove that if $M$
is a diffuse von Neumann algebra, and $X_{1},\ldots ,X_{n}$ is any
family of generators of $M$, then the non-microstates free entropy
dimension of $X_{1},\ldots ,X_{n}$ is at least $1$. These results
are in spirit related to the work of Aagaard \cite{aagaard:dtdimension}.

\begin{acknowledgement*}
This research is based upon work supported by the Clay Mathematics
Institute, a Sloan Foundation Fellowship, and NSF grant DMS-0102332.
\end{acknowledgement*}

\section{Estimates on Free Entropy Dimension.}

\subsection{Partial conjugate variables.}

Let $X_{1},\ldots ,X_{n}\in B(H)$ be self-adjoint variables, which
we assume to be algebraically free. Let $\Omega \in H$ be a trace-vector
for $W^{*}(X_{1},\ldots ,X_{n})$, and $\eta =(\eta _{1},\ldots ,\eta _{n})\in (H\bar{\otimes }H)^{b}$
be given. Let $B$ be an algebra, algebraically free from $X_{1},\ldots ,X_{n}$.

Define the map $\partial ^{\eta }=\partial _{X_{1},\ldots ,X_{n}}^{\eta }B[X_{1},\ldots ,X_{n}]\to H\bar{\otimes }H$
by\begin{eqnarray*}
\partial ^{\eta }(X_{i}) & = & \eta _{i}\\
\partial ^{\eta }(PQ) & = & P\partial ^{\eta }(Q)+\partial ^{\eta }(P)Q\\
\partial ^{\eta }(b) & = & 0,\quad b\in B.
\end{eqnarray*}
The last equation means that $\partial ^{\eta }$ is a derivation
into $H\bar{\otimes }H$ viewed as a bimodule over $B[X_{1},\ldots ,X_{n}]$
in the obvious way.

Note that if $\eta =\mu +\kappa $, then $\partial ^{\eta }=\partial ^{\mu }+\partial ^{\kappa }.$

\begin{defn}
Define the partial conjugate variable $J_{\eta }(X_{1},\ldots ,X_{n}:B)$
to be the unique vector $\xi \in L^{2}(W^{*}(X_{1},\ldots ,X_{n}))$,
so that\[
\langle \xi ,P\rangle =\langle 1\otimes 1,\partial ^{\eta }(P)\rangle _{H\bar{\otimes }H},\qquad \forall P\in B[X_{1},\ldots ,X_{n}],\]
if such a vector exists. If $Y_{1},\ldots ,Y_{m}\in B(H)$ are self-adjoint,
we write $J_{\eta }(X_{1},\ldots ,X_{n}:Y_{1},\ldots ,Y_{m})$ for
$J_{\eta }(X_{1},\ldots ,X_{n}:\mathbb{C}[Y_{1},\ldots ,Y_{n}])$.
\end{defn}
In other words, we set $J_{\eta }(X_{1},\ldots ,X_{n})=(\partial ^{\eta })^{*}(1\otimes 1)$.

Note that, once again, if $\eta =\mu +\kappa $, then\[
J_{\eta }(X_{1},\ldots X_{n}:B)=J_{\mu }(X_{1},\ldots X_{n}:B)+J_{\kappa }(X_{1},\ldots X_{n}:B).\]

\subsection{Examples of partial conjugate variables.}

\begin{example}
\label{example:conjugate}Let $X_{1},\ldots ,X_{n}$ be given. Let
$M=W^{*}(X_{1},\ldots ,X_{n})$, and let $\eta ^{1},\ldots ,\eta ^{n}\in (L^{2}(M)\bar{\otimes }L^{2}(M))^{n}$,
$\eta _{i}^{j}=\sum _{k}a_{k}^{ij}\otimes b_{k}^{ij}$, $a_{k}^{ij},b_{k}^{ij}\in M$
(where the sum is in the $L^{2}$-sense). Let $S_{1},\ldots ,S_{n}$
be a free semicircular family, free from $M$. Let $N=W^{*}(M,S_{1},\ldots ,S_{n})$
and set\[
Y_{j}=\sum _{i=1}^{n}\sum _{k}b_{k}^{ij}S_{i}a_{k}^{ij}\in L^{2}(N).\]
The series defining $Y_{j}$ converges in $L^{2}$ because $S\overline{MS_{j}M}^{L^{2}}\cong L^{2}(M)\bar{\otimes }L^{2}(M)$
in a way that preserves the $L^{2}$ norms. Moreover\[
Y_{j}=J_{\eta ^{j}}(S_{1},\ldots ,S_{n}:X_{1},\ldots ,X_{n}).\]
Indeed, for $w_{j}\in \mathbb{C}[X_{1},\ldots ,X_{n}]$ and $i_{1},\ldots ,i_{n}\in \{1,\ldots ,n\}$\begin{eqnarray*}
\tau (Y_{j}w_{0}S_{i_{1}}w_{1}\cdots w_{k-1}S_{i_{k}}w_{k}) & = & \sum _{i}\sum _{s}\tau (b_{s}^{ij}S_{i}a_{s}^{ij}w_{0}S_{i_{1}}w_{1}\cdots w_{k-1}S_{i_{k}}w_{k})\\
 & = & \sum _{i}\sum _{p,\quad i_{p}=i}\tau (a_{s}^{i_{p}j}w_{0}S_{i_{1}}w_{1}\cdots w_{p-1})\tau (w_{p}S_{i_{p}}\cdots S_{i_{k}}w_{k}b_{s}^{i_{p}j})\\
 & = & \langle 1\otimes 1,\partial _{\eta ^{j}}(w_{0}S_{i_{1}}w_{1}\cdots w_{k-1}S_{i_{k}}w_{k})).
\end{eqnarray*}

\end{example}
\begin{lem}
\label{Lemma:ExampleA}With the above notation, we have\[
J_{\eta ^{j}}(X_{1}+\sqrt{\varepsilon }S_{1},\ldots ,X_{n}+\sqrt{\varepsilon }S_{n})=E_{W^{*}(X_{1}+\sqrt{\varepsilon }S_{1},\ldots ,X_{n}+\sqrt{\varepsilon }S_{n})}(\frac{1}{\sqrt{\varepsilon }}Y_{j}).\]

\end{lem}
\begin{proof}
For all polynomials $P$ in $\mathbb{C}[X_{1}+\sqrt{\varepsilon }X_{1},\ldots ,X_{n}+\sqrt{\varepsilon }S_{n})\subset \mathbb{C}[X_{1},\ldots ,X_{n},S_{1},\ldots ,S_{n}]$,\begin{eqnarray*}
\tau (E_{W^{*}(X_{1}+\sqrt{\varepsilon }S_{1},\ldots ,X_{n}+\sqrt{\varepsilon }S_{n})}(\frac{1}{\sqrt{\varepsilon }}Y_{j})P) & = & \tau (\frac{1}{\sqrt{\varepsilon }}Y_{j}P)\\
 & = & \langle 1\otimes 1,\frac{1}{\sqrt{\varepsilon }}\partial _{(S_{1},\ldots ,S_{n})}^{\eta ^{j}}(P)\rangle \\
 & = & \langle 1\otimes 1,\partial _{(X_{1}+\sqrt{\varepsilon }S_{1},\ldots ,X_{n}+\sqrt{\varepsilon }S_{n})}^{\eta ^{j}}(P)\rangle ,
\end{eqnarray*}
since $\frac{1}{\sqrt{\varepsilon }}\partial ^{\eta ^{j}}(X_{i}+\sqrt{\varepsilon }S_{i})=\eta _{i}^{j}$,
and $\partial $ is a derivation.
\end{proof}
\begin{example}
With the same notation as in Example \ref{example:conjugate}, assume
that $\xi _{j}=J_{\eta ^{j}}(X_{1},\ldots ,X_{n})$ exists. Then\[
\xi _{j}=J_{\eta ^{j}}(X_{1},\ldots ,X_{n}:S_{1},\ldots ,S_{n}).\]

\end{example}
Let $w_{j}\in \mathbb{C}[S_{1},\ldots ,S_{n}]$ and $i_{1},\ldots ,i_{k}\in \{1,\ldots ,n\}$.
We must show that\begin{equation}
\tau (\xi _{j}w_{0}X_{i_{1}}w_{1}\cdots w_{i_{k-1}}X_{i_{k}}w_{k})=\langle 1\otimes 1,\partial _{X_{j}}^{\eta ^{j}}(w_{0}X_{i_{1}}w_{1}\cdots w_{i_{k-1}}X_{i_{k}}w_{k})\rangle .\label{eqn:forConjugate}\end{equation}
The proof is by induction on $k$. If $k=0$, we have\[
\tau (\xi _{j}w_{0})=\tau (\xi _{j})\tau (w_{0}),\]
since $\xi _{j}\in L^{2}(W^{*}(X_{1},\ldots ,X_{n}))$ and $w_{0}$
is free from $W^{*}(X_{1},\ldots ,X_{n})$. Now,\[
\tau (\xi _{j})=\tau (\xi _{j}1)=\langle 1\otimes 1,\partial _{X_{j}}^{\eta ^{j}}(1)\rangle =0.\]
On the other hand, $\partial _{X_{j}}^{\eta ^{j}}(w_{0})=0$, so that
(\ref{eqn:forConjugate}) holds for $k=0$.

Assume that equality in (\ref{eqn:forConjugate}) holds for $k-1$.
By rewriting $W=w_{0}X_{i_{1}}w_{1}\cdots w_{i_{k-1}}X_{i_{k}}w_{k}$
modulo shorter words, we may assume that $W$ has the form\[
W=v_{0}x_{1}v_{1}x_{2}\cdots x_{k}v_{k},\]
where $x_{j}\in \mathbb{C}[X_{1},\ldots ,X_{n}]$, $v_{j}\in \mathbb{C}[S_{1},\ldots ,S_{n}]$,
$\tau (x_{j})=0$ for all $j$, $\tau (v_{j})=0$ for all $j\neq 0$
and $j\neq k$, and either $\tau (v_{0})=0$ or $v_{0}\in \mathbb{C}$,
and either $\tau (v_{k})=0$, or $v_{k}\in \mathbb{C}$. If $k=1$
and $v_{0}=v_{k}\in \mathbb{C}$, the desired equality follows from
the fact that $\xi _{j}=J_{\eta ^{j}}(X_{1},\ldots ,X_{n})$. If not,
$\tau (\xi _{j}W)=0$ by the freeness condition.

On the other hand, we get\[
\langle 1\otimes 1,\partial _{X_{j}}^{\eta }(W)\rangle =\langle 1\otimes 1,\sum _{p}v_{0}x_{1}\cdots v_{p-1}\partial ^{\eta ^{j}}(x_{p})v_{p}\cdots v_{k-1}x_{k}v_{k}\rangle .\]
Let $\partial ^{\eta ^{j}}(x_{p})=\sum _{q}y_{q}^{(1)}\eta _{i_{q}j}y_{q}^{(2)}$,
$y_{q}^{(1)},y_{q}^{(2)}\in \mathbb{C}[X_{1},\ldots ,X_{n}]$. Then
we get, writing $\eta _{i}^{j}=\sum a_{s}\otimes b_{s}$ (series convergent
in $L^{2}$) with $a_{s},b_{s}\in W^{*}(X_{1},\ldots ,X_{n})$:\begin{eqnarray*}
\langle 1\otimes 1,\partial ^{\eta _{j}}(W)\rangle  & = & \sum _{p}\sum _{q}\langle 1\otimes 1,_{p}v_{0}x_{1}\cdots v_{p-1}y_{q}^{(1)}\eta _{i_{q}}^{j}y_{q}^{(2)}v_{p}\cdots v_{k-1}x_{k}v_{k}\rangle \\
 & = & \sum _{p}\sum _{q}\sum _{i}\tau (v_{0}x_{1}\cdots v_{p-1}y_{q}^{(1)}a_{i})\tau (b_{i}y_{q}^{(2)}v_{p}\cdots v_{k-1}x_{k}v_{k}).
\end{eqnarray*}
But each term $\tau (v_{0}x_{1}\cdots v_{p-1}y_{q}^{(1)}a_{i})$ and
$\tau (b_{i}y_{q}^{(2)}v_{p}\cdots v_{k-1}x_{k}v_{k})$ is $0$ by
the freeness condition.

\begin{lem}
\label{lemma:deformation}Let $M=W^{*}(X_{1},\ldots ,X_{n})$, and
let $\eta _{i}^{j}\in L^{2}(M)\bar{\otimes }L^{2}(M)$. Assume that
$\xi _{j}=J_{\eta ^{j}}(X_{1},\ldots ,X_{n})$ exists for all $j=1,\ldots ,n$.
Assume that $S_{1},\ldots ,S_{n}$ are a free semicircular system,
free from $X_{1},\ldots ,X_{n}$. Then $J_{\eta ^{j}}(X_{1}+\sqrt{\varepsilon }S_{1},\ldots ,X_{n}+\sqrt{\varepsilon }S_{n})$
exists and\[
J_{\eta ^{j}}(X_{1}+\sqrt{\varepsilon }S_{1},\ldots ,X_{n}+\sqrt{\varepsilon }S_{n})=E_{W^{*}(X_{1}+\sqrt{\varepsilon }S_{1},\ldots ,X_{n}+\sqrt{\varepsilon }S_{n})}(\xi _{j}).\]

\end{lem}
The proof is similar to that of Lemma \ref{Lemma:ExampleA}.

\begin{prop}
\label{prop:dualSystem}Let $X_{1},\ldots ,X_{n}\in (M,\tau )$, and
let $\Omega \in L^{2}(M)$ be the trace vector. Assume that $D\in B(L^{2}(M))$
is an operator, so that\[
[X_{i},D]=\Xi _{i}\]
and $\Xi _{i}$ belongs to the ideal $HS$ of Hilbert-Schmidt operators
on $L^{2}(M)$. Let $\Psi $ be the identification of $HS$ with $L^{2}(M)\bar{\otimes }L^{2}(M)$
as described in \cite[Proposition 5.11]{dvv:entropy5}, and let $\eta _{i}=\Psi ^{-1}(\Xi _{i})$,
$\eta =(\eta _{1},\ldots ,\eta _{n})$. Then\[
J_{\eta }(X_{1},\ldots ,X_{n})=(D-JD^{*}J)\cdot \Omega .\]

\end{prop}
\begin{proof}
Let $P\in \mathbb{C}[X_{1},\ldots ,X_{n}]$, $P=P^{*}$. Let $\xi =(D-JD^{*}J)\Omega $.
Then\begin{eqnarray*}
\tau (\xi P) & = & \langle P\Omega ,D\Omega \rangle -\langle P\Omega ,JD^{*}J\Omega \rangle =\langle \Omega ,PD\Omega \rangle -\langle D^{*}\Omega ,JP\Omega \rangle \\
 & = & \langle \Omega ,PD\Omega \rangle -\langle D^{*}\Omega ,P\Omega \rangle =\langle \Omega ,PD\Omega \rangle -\langle \Omega ,DP\Omega \rangle \\
 & = & \langle \Omega ,[P,D]\rangle =\Tr (e_{\Omega }[P,D])\\
 & = & \langle 1\otimes 1,\Psi ^{-1}([P,D])\rangle =\langle 1\otimes 1,\partial ^{\eta }(P)\rangle ,
\end{eqnarray*}
because $[\cdot ,D]$ is a derivation, and $\Psi ^{-1}([X_{i},D])=\Psi ^{-1}(\Xi _{i})=\eta _{i}=\partial ^{\eta }(X_{i})$.
\end{proof}

\subsection{Estimates on $\chi ^{*}$ and $\delta ^{*}$.}

We now turn to the task of estimating free Fisher information of $X_{1},\ldots ,X_{n}$,
knowing that $J_{\eta }$ exist. Recall from \cite{dvv:entropy5}
that $\Phi ^{*}(Z_{1},\ldots ,Z_{n})=\sum _{j}\Vert J_{I^{j}}(Z_{1},\ldots ,Z_{n})\Vert _{2}^{2}$,
where $I\in M_{n\times n}(L^{2}(M)\bar{\otimes }L^{2}(M))$ denote
the matrix with entries $I_{i}^{j}=\delta _{ij}1\otimes 1$, and $I^{j}\in (L^{2}(M)\bar{\otimes }L^{2}(M))^{n}$
denotes the row vector $(0,\ldots ,1\otimes 1,\ldots ,0)$, with $1\otimes 1$
in the $j$-th position. For the foregoing, it is convenient to endow
$M_{n\times n}(L^{2}(M)\bar{\otimes }L^{2}(M))$ with the non-normalized
{}``Hilbert-Schmidt'' norm,\[
\Vert T_{ij}\Vert _{M_{n}(L^{2}(M)\bar{\otimes }L^{2}(M))}=\sum _{i,j=1}^{n}\Vert T_{ij}\Vert _{L^{2}(M)\bar{\otimes }L^{2}(M)}.\]

\begin{thm}
\label{thrm:phiestimate}Let $X_{1},\ldots ,X_{n}$ be self-adjoint,
$M=W^{*}(X_{1},\ldots ,X_{n})$ and assume that $\tau $ is a normal
faithful trace on $M$. For each $i,j$, let $\eta _{i}^{j}\in L^{2}(M)\bar{\otimes }L^{2}(M)$
be a vector. Let $\eta ^{j}=(\eta _{1}^{j},\ldots ,\eta _{n}^{j})$.
Assume that for each $j$,\[
\xi _{j}=J_{\eta ^{j}}(X_{1},\ldots ,X_{n})\]
exists. Let $\eta \in M_{n}(L^{2}(M)\bar{\otimes }L^{2}(M))$ be the
matrix with entries $\eta _{i}^{j}$. Let $S_{1},\ldots ,S_{n}$ be
a free semicircular family, free from $M$. 

Then\begin{eqnarray*}
\Phi ^{*}(X_{1}+\sqrt{\varepsilon }S_{1},\ldots ,X_{n}+\sqrt{\varepsilon }S_{n})\leq \sum _{j}\Vert \xi _{j}\Vert _{L^{2}(M)}^{2}+\frac{1}{\varepsilon }\Vert I-\eta \Vert _{M_{n}(L^{2}(M)\bar{\otimes }L^{2}(M))}^{2} &  & \\
+\frac{2}{\sqrt{\varepsilon }}\left(\sum _{j}\Vert \xi _{j}\Vert _{L^{2}(M)}^{2}\right)^{1/2}\Vert I-\eta \Vert _{M_{n}(L^{2}(M)\bar{\otimes }L^{2}(M)).} &  & 
\end{eqnarray*}

\end{thm}
\begin{proof}
Set $X_{j}^{\varepsilon }=X_{j}+\sqrt{\varepsilon }S_{j}$ and $M_{\varepsilon }=W^{*}(X_{1}^{\varepsilon },\ldots ,X_{n}^{\varepsilon })$.
Let $I^{j}\in (L^{2}(M)\bar{\otimes }L^{2}(M))^{n}$ be the vector
$(0,\ldots ,1\otimes 1,\ldots ,0)$. 

From Lemma \ref{Lemma:ExampleA} we get that\[
J_{I^{j}-\eta ^{j}}(X_{1}^{\varepsilon },\ldots ,X_{n}^{\varepsilon })=E_{M_{\varepsilon }}(\frac{1}{\sqrt{\varepsilon }}Y_{j}),\]
with $\Vert Y_{j}\Vert _{2}^{2}=\sum _{i}\Vert \delta _{ij}(1\otimes 1)-\eta _{ij}\Vert _{L^{2}(M)\bar{\otimes }L^{2}(M)}^{2}$.
Hence\[
\sum _{j}\Vert J_{I^{j}-\eta ^{j}}(X_{1}^{\varepsilon },\ldots ,X_{n}^{\varepsilon })\Vert _{L^{2}(M_{\varepsilon })}^{2}\leq \frac{1}{\varepsilon }\Vert I-\eta \Vert _{M_{n}(L^{2}(M)\bar{\otimes }L^{2}(M))}^{2}.\]
From Example \ref{example:conjugate}(b) and Lemma \ref{lemma:deformation}
we get that\begin{eqnarray*}
J_{\eta ^{j}}(X_{1}^{\varepsilon },\ldots ,X_{n}^{\varepsilon }) & = & E_{M_{\varepsilon }}(\xi _{j}),\\
\Vert J_{\eta _{j}}(X_{j}^{\varepsilon }:X_{1}^{\varepsilon },\ldots ,\hat{X_{j}^{\varepsilon }},\ldots ,X_{n}^{\varepsilon })\Vert _{L^{2}(M_{\varepsilon })} & \leq  & \Vert \xi _{j}\Vert _{L^{2}(M)}.
\end{eqnarray*}
Moreover,\[
J_{I^{j}}(X_{1}^{\varepsilon },\ldots ,X_{n}^{\varepsilon })=J_{I^{j}-\eta ^{j}}(X_{1}^{\varepsilon },\ldots ,X_{n}^{\varepsilon })+J_{\eta ^{j}}(X_{1}^{\varepsilon },\ldots ,X_{n}^{\varepsilon }).\]
Let $\Phi =\sum _{j}\Vert J_{I^{j}}(X_{1}^{\varepsilon },\ldots ,X_{n}^{\varepsilon })\Vert _{2}^{2}$.
Then\begin{eqnarray*}
\Phi  & = & \sum _{j}\langle J_{I^{j}}(X_{1}^{\varepsilon },\ldots ,X_{n}^{\varepsilon }),J_{I^{j}}(X_{1}^{\varepsilon },\ldots ,X_{n}^{\varepsilon })\rangle \\
 & = & \sum _{j}\Vert J_{I^{j}-\eta ^{j}}(X_{1}^{\varepsilon },\ldots ,X_{n}^{\varepsilon })\Vert ^{2}+\sum _{j}\Vert J_{\eta ^{j}}(X_{1}^{\varepsilon },\ldots ,X_{n}^{\varepsilon })\Vert ^{2}\\
 &  & +2\sum _{j}\textrm{Re}\langle J_{I^{j}-\eta ^{j}}(X_{1}^{\varepsilon },\ldots ,X_{n}^{\varepsilon }),J_{\eta ^{j}}(X_{1}^{\varepsilon },\ldots ,X_{n}^{\varepsilon })\rangle \\
 & \leq  & \varepsilon ^{-1}\Vert I-\eta \Vert _{M_{n}}^{2}+\sum _{j}\Vert \xi _{j}\Vert ^{2}+2\Vert I-\eta \Vert _{M_{n}}\left[\sum _{j}\Vert J_{\eta ^{j}}(X_{1}^{\varepsilon },\ldots ,X_{n}^{\varepsilon })\Vert ^{2}\right]^{\frac{1}{2}}\\
 & \leq  & \varepsilon ^{-1}\Vert I-\eta \Vert _{M_{n}}^{2}+2\varepsilon ^{-\frac{1}{2}}\Vert I-\eta \Vert _{M_{n}}(\sum _{j}\Vert \xi _{j}\Vert ^{2})^{1/2}+\sum _{j}\Vert \xi _{j}\Vert ^{2},
\end{eqnarray*}
as claimed.
\end{proof}
\begin{cor}
\label{corr:chiAndDeltaEstimate}Under the hypothesis of Theorem \ref{thrm:phiestimate},
there exists a constant $K$ (depending on $n$, $\Vert \xi _{j}\Vert $,
$j=1,\ldots ,n$, and $\Vert X_{j}\Vert _{2}$, $j=1,\ldots ,n$)
so that we have\[
\chi ^{*}(X_{1}+\sqrt{\varepsilon }S_{1},\ldots ,X_{n}+\sqrt{\varepsilon }S_{n})\geq \Vert I-\eta \Vert _{M_{n}(L^{2}(M)\bar{\otimes }L^{2}(M))}^{2}\log \varepsilon ^{1/2}+K,\qquad 0<\varepsilon <\frac{1}{2}.\]
\[
\delta ^{*}(X_{1},\ldots ,X_{n})\geq n-\Vert I-\eta \Vert _{M_{n}(L^{2}(M)\bar{\otimes }L^{2}(M))}^{2}.\]

\end{cor}
\begin{proof}
Let $X_{j}^{\varepsilon }=X_{j}+\sqrt{\varepsilon }S_{j}$. Let also
$\delta =\Vert I-\eta \Vert _{M_{n}}^{2}$, $M=\max \{\Vert \xi _{j}\Vert \}$.
By definition \cite{dvv:entropy5},\[
\chi ^{*}=\chi ^{*}(X_{1}+\sqrt{\varepsilon }S_{1},\ldots ,X_{n}+\sqrt{\varepsilon }S_{n})=C+\frac{1}{2}\int _{0}^{\infty }\left(\frac{n}{1+t}-\Phi ^{*}(X_{1}^{t+\varepsilon },\ldots ,X_{n}^{t+\varepsilon })\right)dt\]
Let $\Phi(t)=\Phi^*(X_1^t,\ldots,X_n^t)$.
We have\begin{eqnarray*}
\chi ^{*} & = & C+\frac{1}{2}\int _{0}^{\infty }\left[\frac{n}{1+t}-\Phi({t+\varepsilon })\right]dt\\
 & = & C+\frac{1}{2}\int _{\varepsilon }^{\frac{1}{2}}\left[\frac{n}{1+t-\varepsilon }-\Phi({t})\right]dt+\frac{1}{2}\int _{\frac{1}{2}}^{\infty }\left[\frac{n}{1+t-\varepsilon }-\Phi ({t})\right]dt\\
 & \geq  & C+\frac{1}{2}\int _{\varepsilon }^{\frac{1}{2}}\left[\frac{n}{1+t}-\frac{\delta }{t}-\frac{2\delta ^{\frac{1}{2}}M}{\sqrt{t}}-M^{2}\right]dt+\frac{1}{2}\int _{\frac{1}{2}}^{\infty }\left[\frac{n}{1+t-\varepsilon }-\frac{n}{1+t}\right]dt\\
 & \geq  & C'-n\log (1+\varepsilon )^{1/2}+\delta \log t^{1/2}+4\sqrt{\varepsilon }\delta ^{\frac{1}{2}}M+\varepsilon M^{2}+\frac{n}{2}\log \frac{1+t}{1+t-\varepsilon }\Big |_{\frac{1}{2}}^{\infty }\\
 & \geq  & K+\delta \log t^{1/2}
\end{eqnarray*}
Now, by definition we have\begin{eqnarray*}
\delta ^{*}(X_{1},\ldots ,X_{n}) & = & n+\limsup _{t\to 0}\frac{\chi ^{*}(X_{1}^{\varepsilon },\ldots ,X_{n}^{\varepsilon })}{|\log \varepsilon ^{1/2}|}\\
 & \geq  & n+\limsup _{t\to 0}\frac{K+\delta \log t^{1/2}}{|\log t^{1/2}|}=n-\delta .
\end{eqnarray*}
It follows that $\delta ^{*}(X_{1},\ldots ,X_{n})\geq n-\Vert I-\eta \Vert _{M_{n}}^{2}$.
\end{proof}
\begin{lem}
\label{Lemma:distAndDim}Let $N$ be a finite von Neumann algebra
with a faithful normal trace $\tau $. Let $n$ be a finite integer,
and let $H=L^{2}(N,\tau )^{n}$ be a left module over $N$. Denote
by $\Omega \in L^{2}(N,\tau )$ the GNS vector associated to $\tau $. 

Let $K\subset H$ be a closed $N$-invariant subspace of $H$. Endow
$M_{n\times n}(L^{2}(M))$ with the norm\[
\Vert h_{ij}\Vert _{M_{n}}^{2}=\sum _{ij=1}^{n}\Vert h_{ij}\Vert ^{2}.\]
Let $A(K)=\{T\in M_{n\times n}(M):TH\subset K\}\cong K^{n}$. Then
we have:\[
\dim _{N}K=n-\dist (I,A(K))^{2},\]
where $I\in M_{n}(H)$ denotes the matrix $I_{ij}=\delta _{ij}\Omega $,
and the distance is computed with respect to the Hilbert space norm
on $M_{n\times n}(L^{2}(M))$ described above.
\end{lem}
\begin{proof}
The commutant $N'$ of $N$ acting on $H$ can be identified with
the algebra of $n\times n$ matrices $M_{n}(N)$. Endow this algebra
with the non-normalized trace $\Tr $, defined by the property that
$\Tr (I)=n$, where $I\in M_{n}(N)$ denote the identity matrix. Let
$e_{K}\in N'$ be the orthogonal projection from $H$ onto $K$. Then
$\dim _{N}K=\Tr (e_{K})$. Thus\[
n-\dim _{N}K=\dim _{N}K^{\perp }=\Tr (I-e_{K}).\]
Since $I-e_{K}$ is a projection, we obtain that\[
n-\dim _{N}K=\Tr ((I-e_{K})^{2})=\Vert I-e_{K}\Vert _{L^{2}(M_{n}(N),\Tr )}^{2}.\]
Now, $L^{2}(M_{n}(N),\Tr )=M_{n}(H)$ isometrically. Moreover, since
the orthogonal projection of $I$ onto $A(K)$ is exactly $e_{K}$,
we get that $\Vert I-e_{K}\Vert _{L^{2}(M_{n}(N),\Tr )}$ is exactly
the distance from $I$ to the subspace $M_{n}(K)$.
\end{proof}
\begin{lem}
\label{lemma:bimodule}Let $X_{1},\ldots ,X_{n}\in (M,\tau )$, $M=W^{*}(X_{1},\ldots ,X_{n})$.
Let\[
H(X_{1},\ldots ,X_{n})=\{\eta =(\eta _{1},\ldots ,\eta _{n})\in L^{2}(M,\tau )^{n}:J_{\eta }(X_{1},\ldots ,X_{n})\textrm{ exists}\}.\]
Then $H(X_{1},\ldots ,X_{n})$ is a bimodule over $J\mathbb{C}[X_{1},\ldots ,X_{n}]J$
and its closure is a $JMJ\bar{\otimes }MJ$-submodule of $(L^{2}(M)\bar{\otimes }L^{2}(M))^{n}$. 
\end{lem}
\begin{proof}
Clearly, only the first assertion needs to be proved. 

Let $\eta =(\eta _{1},\ldots ,\eta _{n})\in H(X_{1},\ldots ,X_{n})$.
Thus $\xi =J_{\eta }(X_{1},\ldots ,X_{n})=\partial _{\eta }^{*}(1\otimes 1)$
exists. Let $a,b,P\in \mathbb{C}[X_{1},\ldots ,X_{n}]$. Then\begin{eqnarray*}
\langle 1\otimes 1,\partial ^{Jb^{*}J\eta Ja^{*}J}(P)\rangle  & = & \langle a\otimes b,\partial ^{\eta }(P)\rangle .
\end{eqnarray*}
We now have, exactly as in the proof of Proposition 4.1 in \cite{dvv:entropy5},\[
\langle 1\otimes a,\partial ^{\eta }(P)\rangle =\langle \xi a-(\id \otimes \tau )((\partial ^{\eta }(a^{*})^{*}),P\rangle .\]
Hence\[
J_{\eta Ja^{*}J}(X_{1},\ldots ,X_{n})=J_{\eta }(X_{1},\ldots ,X_{n})a-(\id \otimes \tau )((\partial ^{\eta }(a^{*}))^{*}.\]
Similarly, we have $J_{Jb^{*}J\eta }=bJ_{\eta }(X_{1},\ldots ,X_{n})-(\tau \otimes \id )((\partial ^{\eta }(a^{*}))^{*}$.
\end{proof}
\begin{cor}
\label{corr:Hestimate}Let $X_{1},\ldots ,X_{n}\in (M,\tau )$, $M=W^{*}(X_{1},\ldots ,X_{n})$.
With the notation of Lemma \ref{lemma:bimodule}, we have\[
\delta ^{*}(X_{1},\ldots ,X_{n})\geq \dim _{JMJ\bar{\otimes }JMJ}\overline{H(X_{1},\ldots ,X_{n})}.\]
Moreover, if $\eta ^{j}=(\eta _{1}^{j},\ldots ,\eta _{n}^{j})\in H(X_{1},\ldots ,X_{n})$
for all $j=1,\ldots ,n$, then we have\[
\dim _{JMJ\bar{\otimes }JMJ}\overline{H(X_{1},\ldots ,X_{n})}\geq n-\Vert I-(\eta _{i}^{j})\Vert _{M_{n}(L^{2}(M)\bar{\otimes }L^{2}(M))}^{2}.\]

\end{cor}
\begin{proof}
We apply Lemma \ref{Lemma:distAndDim} with $N=JMJ\bar{\otimes }JMJ$,
$H=(L^{2}(M)\bar{\otimes }L^{2}(M))^{n}$, and $K=\overline{H(X_{1},\ldots ,X_{n})}$.
The space $A(K)$ can then be identified with the closure of the set\[
A=\{(\eta _{i}^{j})\in M_{n}(L^{2}(M)\bar{\otimes }L^{2}(M)):(\eta _{1}^{j},\ldots ,\eta _{n}^{j})\in H(X_{1},\ldots ,X_{n})\  \forall j\}.\]
Applying Corollary \ref{corr:chiAndDeltaEstimate}, we find that\[
\delta ^{*}(X_{1},\ldots ,X_{n})\geq n-\dist (I,A)^{2}.\]
By Lemma \ref{Lemma:distAndDim}, we finally get that \[
\delta ^{*}(X_{1},\ldots ,X_{n})\geq \dim _{N}K,\]
as claimed. 

The last inequality follows from the fact that $(\eta _{i}^{j})\in A$
and hence one has $\Vert I-(\eta _{i}^{j})\Vert _{M_{n}(L^{2}(M)\bar{\otimes }L^{2}(M))}^{2}\geq \dist (I,A)^{2}$.
\end{proof}
\begin{cor}
\label{corr:HSdeltaEstimate}Let $X_{1},\ldots ,X_{n}\in (M,\tau )$,
$M=W^{*}(X_{1},\ldots ,X_{n})$. Let $HS$ be the space of Hilbert-Schmidt
operators on $L^{2}(M,\tau )$. Set\[
H_{0}(X_{1},\ldots ,X_{n})=\{(\Xi _{1},\ldots ,\Xi _{n})\in HS^{n}:\exists D\in B(L^{2}(M,\tau ))\textrm{ s.t. }\Xi _{i}=[D,X_{i}]\  \forall i\}.\]
Then\[
\delta ^{*}(X_{1},\ldots ,X_{n})\geq \dim _{JMJ\bar{\otimes }JMJ}\overline{H_{0}(X_{1},\ldots ,X_{n})}.\]
Furthermore, if $D_{1},\ldots ,D_{n}\in B(L^{2}(M,\tau ))$ are such
that $\Xi _{ij}=[D_{i},X_{j}]\in HS$, then\[
\dim _{JMJ\bar{\otimes }JMJ}\overline{H_{0}(X_{1},\ldots ,X_{n})}\geq n-\Vert I-(\Xi _{ij})\Vert _{M_{n}(L^{2}(M)\bar{\otimes }L^{2}(M))}^{2}.\]

\end{cor}
\begin{proof}
By Proposition \ref{prop:dualSystem}, $H_{0}(X_{1},\ldots ,X_{n})\subset H(X_{1},\ldots ,X_{n})$.
Note that if $\Xi _{j}=[D,X_{j}]$, $j=1,\ldots ,n$, then $JaJ\Xi _{j}JbJ=[JaJDJbJ,X_{j}]$
for any $a,b\in M$. Thus $H_{0}(X_{1},\ldots ,X_{n})$ is an $M,M$-bimodule.
The rest follows from Corollary \ref{corr:Hestimate}.
\end{proof}
This fact has the following amusing consequence, which is the non-microstates
analog of the {}``hyperfinite monotonicity'' of \cite{kenley:hyperfinite}.

\begin{thm}
Assume that $X_{1},\ldots ,X_{n}$ are self-adjoint and generate a
von Neumann algebra $M$ with a normal faithful trace $\tau $. Assume
that $M$ is diffuse. Then $\delta ^{*}(X_{1},\ldots ,X_{n})\geq 1$.
\end{thm}
\begin{proof}
We note that the map\[
HS\ni D\mapsto ([D,X_{1}],\ldots ,[D,X_{n}])\in H_{0}(X_{1},\ldots ,X_{n})\]
is injective. Since $\dim _{JMJ\bar{\otimes }JMJ}HS=1$, we find that
$\dim _{JMJ\bar{\otimes }JMJ}\overline{H_{0}(X_{1},\ldots ,X_{n})}\geq 1$.
Applying Corollary \ref{corr:HSdeltaEstimate}, we get that $\delta ^{*}(X_{1},\ldots ,X_{n})\geq 1$.
\end{proof}
We note that the estimate in Corollary \ref{corr:HSdeltaEstimate}
is optimal for a single variable $X$. The map\[
HS(L^{2}(X))\ni D\mapsto [X,D]\]
has a kernel exactly if the spectral measure of $X$ has atoms. Moreover,
if we denote by $\mu $ the distribution of $X$, the Murray-von Neumann
dimension over $W^{*}(X)\bar{\otimes }W^{*}(X)$ of its kernel is
given by $\sum _{t\in \mathbb{R}}\mu (\{t\})^{2}.$ From this we get
that $\dim _{W^{*}(X)\bar{\otimes }W^{*}(X)}\overline{H_{0}(X)}=1-\sum _{t\in \mathbb{R}}\mu (\{t\})^{2}$.
Since for a single variable, $\delta =\delta ^{*}$, we see that this
estimate is optimal.

\section{Free dimension for $q$-Semicircular variables.}

We recall the construction of $q$-semicircular variables given by
Bozejko and Speicher in \cite{Speicher:q-comm}

Let $N$ be an integer, $H=\mathbb{C}^{N}$, and $-1<q<1$. Consider
the vector space\[
F_{alg}(H)=\mathbb{C}\Omega \oplus \bigoplus _{n\geq 1}H^{\otimes n}\]
(algebraic direct sum). This vector space is endowed with a positive
definite inner product $\langle \cdot ,\cdot \rangle $, which depends
on $q$, and which we do not describe here. Denote by $F_{q}(H)$
the completion of $F_{alg}(H)$ with respect to this inner product.

For $h\in H$, define $\ell (h):F_{q}(H)\to F_{q}(H)$ by extending
continuously the map\begin{eqnarray*}
\ell (h)h_{1}\otimes \cdots \otimes h_{n} & = & h\otimes h_{1}\otimes \cdots \otimes h_{n},\\
\ell (h)\Omega  & = & h.
\end{eqnarray*}
The adjoint is given by\begin{eqnarray*}
\ell ^{*}(h)h_{1}\otimes \cdots \otimes h_{n} & = & \sum _{k=1}^{n}q^{k-1}\langle h_{k},h\rangle h_{1}\otimes \cdots \otimes \hat{h}_{k}\otimes \cdots \otimes h_{n},\\
\ell ^{*}(h)\Omega  & = & 0,
\end{eqnarray*}
where $\hat{\cdot }$ denotes omission.

Consider also $r(h)$ given by\begin{eqnarray*}
r(h)h_{1}\otimes \cdots \otimes h_{n} & = & h_{1}\otimes \cdots \otimes h_{n}\otimes h\\
r(h)\Omega  & = & h.
\end{eqnarray*}
Finally, let $P_{n}:F_{q}(H)\to F_{q}(H)$ be the orthogonal projection
onto tensors of rank $n$. Let $\Xi _{q}=\sum _{n\geq 0}q^{n}P_{n}$.

\begin{lem}
\label{lemma:commutators}Let $|q|<1$. With the above notation we
have\begin{eqnarray*}
{}[\ell (h),r(g)] & = & 0,\\
{}[\ell (h)^{*},r(g)] & = & \langle g,h\rangle \Xi _{q}.
\end{eqnarray*}
In particular, $[\ell (h)^{*},r(g)]$ is compact if $N=\dim H<\infty $.
Moreover, if $q^{2}N<1$, then $\Xi _{q}$ is Hilbert-Schmidt, and
$\Vert P_{0}-\Xi _{q}\Vert _{HS}^{2}=\frac{q^{2}N}{1-q^{2}N}.$
\end{lem}
\begin{proof}
We have\[
[\ell (h),r(g)]\Omega =h\otimes g-h\otimes g=0.\]
Similarly,\[
[l(h),r(g)]h_{1}\otimes \cdots \otimes h_{n}=h\otimes h_{1}\otimes \cdots \otimes h_{n}\otimes g-h\otimes h_{1}\otimes \cdots \otimes h_{n}\otimes g=0.\]
Also,\[
[\ell (h)^{*},r(g)]\Omega =\ell ^{*}(h)g-0=\langle g,h\rangle \Omega =\langle g,h\rangle P_{0}\Omega ,\]
and for $\xi =h_{1}\otimes \cdots \otimes h_{n}$,\begin{eqnarray*}
[\ell (h)^{*},r(g)]\xi  & = & \ell (h)^{*}(\xi \otimes g)-r(g)\sum _{k=1}^{n}q^{k-1}\langle h_{k},h\rangle h_{1}\otimes \cdots \otimes \hat{h}_{k}\otimes \cdots \otimes h_{n}\\
 & = & \sum _{k=1}^{n}q^{k-1}\langle h_{k},h\rangle h_{1}\otimes \cdots \otimes \hat{h}_{k}\otimes \cdots \otimes h_{n}\otimes g+q^{n}\langle g,h\rangle \xi \\
 & - & \sum _{k=1}^{n}q^{k-1}\langle h_{k},h\rangle h_{1}\otimes \cdots \otimes \hat{h}_{k}\otimes \cdots \otimes h_{n}\otimes g\\
 & = & q^{n}\langle g,h\rangle \xi =q^{n}\langle g,h\rangle P_{n}(\xi ).
\end{eqnarray*}
Since the rank of $P_{n}$ is finite if $\dim H<\infty $, and $q^{n}\to 0$
as $n\to \infty $, we find that $\Xi _{q}=\sum q^{n}P_{n}$ is compact.

Since $P_{n}$ are orthogonal for the Hilbert-Schmidt norm, we find
that the square of the Hilbert-Schmidt norm of $P_{0}-\Xi _{q}$ is
given by\begin{eqnarray*}
\Vert P_{0}-\Xi _{q}\Vert _{HS}^{2} & = & \sum _{n\geq 1}q^{2n}\Tr (P_{n})\\
 & = & \sum _{n\geq 1}(q^{2}N)^{n}=\frac{q^{2}N}{1-q^{2}N},
\end{eqnarray*}
if $q^{2}N<1$, and is infinite otherwise.
\end{proof}
Let $e_{1},\ldots ,e_{n}$ be an orthonormal basis for $\mathbb{C}^{N}$.
Applying Corollary \ref{corr:HSdeltaEstimate}, to $X_{1}=\ell (e_{1})+\ell (e_{1})^{*},\ldots ,X_{N}=\ell (e_{N})+\ell (e_{N})^{*}$,
$D_{1}=r(e_{1}),\ldots ,D_{N}=r(e_{N})$, $\Xi _{ij}=\delta _{ij}\Xi _{q}$,
we get:

\begin{cor}
Let $q$ and $N$ be so that $q^{2}N<1$. Then\[
\delta ^{*}(X_{1},\ldots ,X_{N})\geq N-\frac{q^{2}N^{2}}{1-q^{2}N}=N\left(1-\frac{q^{2}N}{1-q^{2}N}\right).\]
Moreover, $\delta ^{*}(X_{1},\ldots ,X_{N})>1$ for $q$ so that $q^{2}N<1$.
\end{cor}
The fact that $\delta ^{*}(X_{1},\ldots ,X_{N})>1$ follows from the
fact that $H_{0}(X_{1},\ldots ,X_{n})$ is clearly strictly larger
than the closure of $\{([D,X_{1}],\ldots ,[D,X_{n}]):D\in HS\}$,
since it contains vectors of the form $(\Xi _{q},0,\ldots ,0)$.

We should point out that the estimate above behaves badly for large
values of $q$ and there is no reason to believe that it is optimal,
even for $q$ small. It would be interesting to see if {}``better''
choices for $D_{j}$ can improve this estimate, or to compute the
Murray-von Neumann dimension in Corollary \ref{corr:HSdeltaEstimate}.

We should also point out that the natural question is now whether
$\delta (X_{1},\ldots ,X_{n})$ (here we mean the microstates free
entropy dimension as introduced in \cite{dvv:entropy2}) is bigger
than $1$ for $q$-semicircular families. In conjunction with this,
we mention the results of Sniady \cite{sniady:q-RM-model}, who showed
that $W^{*}(X_{1},\ldots ,X_{n})$ embeds into the ultrapower of the
hyperfinite II$_{1}$ factor (i.e., $X_{1},\ldots ,X_{n}$ have matricial
microstates).

\section{Ozawa's condition and $q$-semicircular random variables.}

In his remarkable paper \cite{ozawa:solid}, Ozawa discusses the consequences
of the following condition for a finite von Neumann algebra $M$.
Let $K$ denote the ideal of compact operators inside $B(L^{2}(M))$,
and denote by $\pi :C^{*}(M,M')\to B(L^{2}(M))/K$ the restriction
of the quotient map.

\begin{condition}
(Ozawa) There are two $C^{*}$-subalgebras $A$ of $M$ and $B$ of
$M'$, so that\\
(a) $A$ is locally reflexive\\
(b) $A$ generates $M$ and $B$ generates $M'$\\
(c) The map\[
\sum a_{i}\otimes b_{i}\mapsto \pi \left(\sum a_{i}b_{i}\right)\in B(L^{2}(M))/K\]
extends to a continuous map on the minimal (spatial) tensor product
$A\otimes _{\min }B$.
\end{condition}
Ozawa proves that if $M$ satisfies this condition, then it is {}``solid'':
the relative commutant $N'\cap M$ of any diffuse von Neumann subalgebra
of $M$ is hyperfinite. Note that in particular, if the center of
$Z(M)$ is diffuse, then it follows that $M$ must be hyperfinite,
since $M=Z(M)'\cap M$. Ozawa also proved that if $M$ is a factor,
and $M$ is non-hyperfinite, then $M$ does not have property $\Gamma $.

\begin{thm}
Let $|q|<\sqrt{2}-1$ and let $X_{1},\ldots ,X_{n}$, $n\geq 2$,
be a family of $q$-semicircular random variables%
\footnote{We could actually use here a better bound given by Dykema and Nica,
of $|q|$ less than approximately $0.44$, see part $4$ of the Theorem
on page 203 of \cite{dykema-nica:stability}.%
}. Then $M=W^{*}(X_{1},\ldots ,X_{n})$ satisfies Ozawa's condition.
In particular:\\
(a) $M$ is solid, i.e., $N'\cap M$ is hyperfinite for any diffuse
subalgebra $N\subset M$;\\
(b) $M\not \cong N_{1}\otimes N_{2}$ with $N_{1}$ and $N_{2}$ $II_{1}$
factors. The center $Z(M)$ is not diffuse.

In addition, if $M$ is a factor, it follows that $M$ is non-$\Gamma $.
\end{thm}
\begin{proof}
Denote by $e_{1},\ldots ,e_{n}$ an orthonormal basis for $\mathbb{C}^{n}$,
so that $X_{j}=\ell (e_{j})+\ell (e_{j})^{*}.$ Let\[
Y_{j}=r(e_{j})+r(e_{j})^{*}.\]
Then $M'=W^{*}(Y_{1},\ldots ,Y_{n})$. Let\[
A=C^{*}(X_{1},\ldots ,X_{n}),\quad B=C^{*}(Y_{1},\ldots ,Y_{n}),\quad C=C^{*}(r(e_{1}),\ldots ,r(e_{n})).\]

By \cite{dykema-nica:stability} and \cite{jorgensen-at-all:stabilityCuntz},
$C$ is isomorphic to an extension of the Cuntz algebra, and is therefore
nuclear. This implies that $A$ (and $B\cong A$) is locally reflexive
(see the discussion in \cite{ozawa:solid}).

Furthermore, note that $[A,C]\subset K$. This follows from Lemma
\ref{lemma:commutators}, where we note that $[\ell (e_{i}),r(e_{j})]=0$
and $[\ell (e_{i}),r(e_{j})^{*}]\in K$.

Denote by $\pi :B(L^{2}(M))\to B(L^{2}(M))/K$ the quotient map. It
follows that the images $\pi (A)$ and $\pi (C)$ commute. Thus the
map\[
\sum a_{i}\otimes b_{i}\mapsto \pi (\sum a_{i}b_{i}),\qquad a_{i}\in A,\  b_{i}\in C\]
is continuous for the max tensor product norm on $A\otimes C$. Since
$C$ is nuclear, we find that the min and max tensor norms on $A\otimes C$
coincide. Thus this map is $\otimes _{\min }$-continuous.

Restricting this map to the image of $A\otimes B\subset A\otimes _{\min }B\subset A\otimes _{\min }C$
we obtain finally that Ozawa's condition is indeed satisfied. Thus
(a) holds. Since by \cite{nou:qNonInjective}, $M$ is not hyperfinite,
(b) also holds.
\end{proof}
Note that the only reason for imposing any restrictions on the values
of $q$ in the argument above is that we need to know that $C^{*}(r(e_{1}),\ldots ,r(e_{n}))$
is nuclear. It may be possible to argue that this is the case directly,
without relying on the stability results \cite{dykema-nica:stability}
and \cite{jorgensen-at-all:stabilityCuntz}. 

Lastly, we should mention that we did not settle the issue of whether
$M$ is a factor, which is only known for certain $q$ \cite{sniady:qfactor}
and also in the case of an infinite family of generators \cite{krolak:q-Gaussian-Factors}.
In the case that these algebras are factors and the number of generators
is finite, the results of Ozawa and Popa \cite{ozawa-popa:primeFactorization}
on tensor powers of $M_{q}$ become available.

\bibliographystyle{amsplain}

\providecommand{\bysame}{\leavevmode\hbox to3em{\hrulefill}\thinspace}

\end{document}